\documentclass{article} 
\usepackage{graphicx}
\usepackage{amsmath}
\usepackage{amsfonts}
\usepackage{amssymb}
\usepackage{latexsym}
\usepackage{tikz}
\usepackage{fancyvrb}
\usepackage{datetime}

\usepackage{graphicx} 
\usepackage{color}
\usepackage[colorlinks]{hyperref}

\newtheorem{lemma}{Lemma}[section]
\newtheorem{proposition}{Proposition}[section]
\newtheorem{corollary}{Corollary}[section]
\newtheorem{remark}{Remark}
\newtheorem{theorem}{Theorem}[section]
\newtheorem{conjecture}{Conjecture}

\usepackage{graphicx}
\usepackage{amsmath}
\usepackage{amsfonts}
\usepackage{amssymb}
\usepackage{latexsym}

\begin{document}
	
	\begin{center}
{\Large{On the spectral determinations of the connected multicone graphs $ K_r\bigtriangledown sK_t $}}

\medskip	

{Ali Zeydi Abdian\footnote{Lorestan University, College of Science, Lorestan, Khoramabad, Iran; e-mail: abdian.al@fs.lu.ac.ir;  aabdian67@gmail.com; azeydiabdi@gmail.com}, Lowell W. Beineke\footnote{Indiana University -- Purdue University Fort Wayne, Fort Wayne, Indiana, U.S.A.; email: beineke@ipfw.edu} and Afshin Behmaram
\footnote{Faculty of Mathematical Sciences, University of Tabriz, Tabriz, Iran; Email: behmaram@tabrizu.ac.ir}}
	
	\end{center}
	
	\begin{abstract}
			
			 In this study we investigate the spectra of the family of connected multicone graphs. A multicone graph is defined to be the join of a clique and a regular graph. Let $ r $, $ t $ and $ s $ be natural numbers, and let $ K_r $  denote a complete graph on $ r $ vertices. It is proved that connected multicone graphs $ K_r\bigtriangledown sK_t $, a natural generalization of friendship graphs, are determined by their adjacency spectra as well as their Laplacian spectra. Also, we show that the complement of multicone graphs $ K_r\bigtriangledown sK_t $ are determined by their adjacency spectra, where $ s\neq 2 $.
	\end{abstract}
	
	\section{ Introduction}
	A long-standing question connecting graph theory and linear algebra has to do with the set of eigenvalues of the adjacency matrix of a graph, called the spectrum of the graph.  Although it is well known that different graphs can have the same spectrum, it remains an open question as to whether most graphs have a spectrum shared by another graph or not.  In fact, not many families of graphs are known that have their own spectrum, not shared by any other graphs. In the past decades, graphs that are determined by their spectrum have received much more and more attention, since they have been applied to several fields, such as randomized algorithms, combinatorial optimization problems and machine learning. An important part of spectral graph theory is devoted to determining whether given graphs or classes of graphs are determined by their spectra or not. So, finding and introducing any class of graphs which are determined by their spectra can be an interesting and important problem.  We begin with some of the notation and terminology that will be used in the paper.  All graphs considered here are simple and undirected, and, in general, given a graph $G$, $n$ will denote the number of vertices (also called its \textit{order}) and $m$ the number of edges. If its vertices are $v_1, v_2, \ldots, v_n$, then its \textit{adjacency matrix} $\mathbf{A}(G)$	is the $n \times n$	matrix with $a_{ij} = 1$ if $v_i$ and $v_j$ are adjacent and $0$ otherwise.  Its \textit{degree matrix} is defined to be the diagonal matrix $\mathbf{D}(G) = {\rm{diag}}(d_1, d_2, \ldots, d_n)$, where $d_i$ is the degree of vertex $v_i$.  Two other matrices are defined in terms of these: the \textit{Laplacian matrix} is $\mathbf{L}(G) = \mathbf{D}(G) - \mathbf{A}(G)$  and the \textit{signless Laplacian matrix} is $\mathbf{S}(G) = \mathbf{D}(G) + \mathbf{A}(G)$.
			We denote the characteristic polynomial $\det (x \mathbf{I} - \mathbf{A})$ of $ G $ by $P_G(x)$. A number $\lambda$ is an \textit{eigenvalue} of $G$ if it is a root of this polynomial.  Since $\mathbf{A}(G)$ is a symmetric matrix, all of its eigenvalues are real.  The \textit{adjacency spectrum} (\textit{Laplacian spectrum}) of $ G $, denoted  $\operatorname{Spec}_{A}(G)$  ($\operatorname{Spec}_{L}(G)$) , is the multiset of these eigenvalues. Two graphs $G$ and $H$ are said to be \textit{ $A$-cospectral ($L$-cospectral)} if the corresponding adjacency spectra (Laplacian spectra) are the same. A graph $G$ is said to be \textit{ $DAS$ ($DLS$)} if there is no other non-isomorphic graph $A$-cospectral ($L$-cospectral) with it, i.e., $\operatorname{Spec}_{A}(H) = \operatorname{Spec}_{A}(G)$ ($\operatorname{Spec}_{L}(H) = \operatorname{Spec}_{L}(G)$) implies $G\cong H$. By analogy, we define \textit{determined by signless Laplacian spectrum} \textit{($DQS$ for short)} graphs. The key question that we consider is the extent to which the spectrum (of either type) of a graph is unique; that is, whether there is only one graph with that spectrum.\\

	So far numerous examples of cospectral but non-isomorphic
graphs have been constructed by interesting techniques such as Seidel switching, Godsil-McKay switching, Sunada or Schwenk method. For more information, one may see \cite{BH, VH, VAH} and the references cited in them. Only a few graphs with very special structures have been reported to be determined by their spectra (DS, for short) (see \cite{BJ, CHVW, DH, HLZ, LS, SZ, W2, WBHB} and the references cited in them). Recently Wei Wang and Cheng-Xian Xu have developed a new method in \cite{W2} to show that many graphs are determined by their spectrum and the spectrum of their complement.\\

			One of the first investigations into this question was made in 1971 by Harary, et al. \cite{HKMR}.  They asserted that (stated in slightly different terminology), based on the data they computed for graphs with up to seven vertices, ``one is tempted to conjecture'' that the fraction of graphs with spectra that are not unique decreases as the order increases.  Technically, this is not exactly the same as the conjecture that the probability goes to $0$, but the two are closely related:\\   
			
			\noindent \textbf {Unique Spectrum Conjecture}  \textit{Almost all graphs are determined by their spectrum.}\\
			
			One fact that makes this conjecture especially intriguing is that there is one very interesting family of graphs for which the corresponding statement is known not to hold.  In fact, Schwenk \cite{Sch} proved that it is about as far off as it could be.\\
			
			\noindent \textbf {Co-spectral Tree Theorem}  \textit{Almost no trees are determined by their spectrum.}\\  
			
			What this means is that, as $n \rightarrow \infty$, the fraction of trees of order $n$ that have the same spectrum as another tree approaches $1$.  
			
			There are of course many versions of a conjecture such as the one above, not only for the different types of spectra, but also for different families of graphs.\\ 
			
			The general terminology that we use may be found in standard textbooks on graph theory, but we give some that will be used here, some of which varies from author to author.  In particular, we use the following notation on graph operations.  We define the \textit{sum} $G + H$ of two vertex-disjoint graphs $G$ and $H$ to be their union; that is, $V(G + H) = V(G) \cup V(H))$ and $E(G + H) = E(G) \cup E(H))$.  Clearly, this can be extended to more graphs, $G_1 + G_2 + \ldots + G_k$, and the sum of $k$ copies of the same graph $G$ is denoted $kG$. The \textit{join} $G \triangledown H$ (or $G \ast H$) is obtained from $G + H$ by adding an edge from each vertex of $G$ to each vertex of $H$, that is, by adding the set of edges $\{vw: v \in V(G), w \in V(H)\}$.\\
			
			The graphs that we consider here are combinations of sums and joins.  We begin with a special case known as a \textit{friendship graph} (also known as a (\textit{Dutch}) \textit{windmill}). Erd\"os, R\'enyi, and S\'os \cite{ERS} proved that if $G$ is the graph of $n$ people for which each pair have exactly one friend in common, then $G$ consists of $t$ triangles (with $n$ odd and $t = \frac{1}{2}(n - 1)$), all having one common vertex. This graph is denoted $F_t$, and $F_2$, $F_3$, and $F_4$ are shown in Fig. 1.

			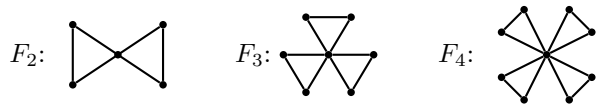
\begin{figure}
				\begin{center}
					\begin{tikzpicture}
					
					\node  at (0,.4)  [scale=1]{$F_2$:};
					
					\node (2a) at (.6,0)  [circle,fill=black,scale=0.3]{};
					\node (2b) at (.6,.8)  [circle,fill=black,scale=0.3]{};
					\node (2c) at (1.2,.4)  [circle,fill=black,scale=0.3]{};
					\node (2d) at (1.8,0)  [circle,fill=black,scale=0.3]{};
					\node (2e) at (1.8,.8)  [circle,fill=black,scale=0.3]{};
					
					\draw [thick] (2a) -- (2e);
					\draw [thick] (2a) -- (2b);
					\draw [thick] (2b) -- (2d);
					\draw [thick] (2d) -- (2e);
					
					\node  at (3,.4)  [scale=1]{$F_3$:};
					
					\node (3a) at (3.7,.9)  [circle,fill=black,scale=0.3]{};
					\node (3b) at (4.3,.9)  [circle,fill=black,scale=0.3]{};
					\node (3c) at (4,.4)  [circle,fill=black,scale=0.3]{};
					\node (3d) at (3.4,.4)  [circle,fill=black,scale=0.3]{};
					\node (3e) at (4.6,.4)  [circle,fill=black,scale=0.3]{};
					\node (3f) at (3.7,-.1)  [circle,fill=black,scale=0.3]{};
					\node (3g) at (4.3,-.1)  [circle,fill=black,scale=0.3]{};
					
					\draw [thick] (3a) -- (3b) -- (3c) -- (3a) ;
					\draw [thick] (3d) -- (3f) -- (3c) -- (3d) ;
					\draw [thick] (3e) -- (3g) -- (3c) -- (3e) ;
					
					\node  at (5.7,.4)  [scale=1]{$F_4$:};
					
					\node (4a) at (6.9, .4)  [circle,fill=black,scale=0.3]{};
					\node (4b) at (6.6,1)  [circle,fill=black,scale=0.3]{};
					\node (4c) at (7.2,1)  [circle,fill=black,scale=0.3]{};
					\node (4d) at (6.3,.7)  [circle,fill=black,scale=0.3]{};
					\node (4e) at (6.3,.1)  [circle,fill=black,scale=0.3]{};
					\node (4f) at (6.6,-.2)  [circle,fill=black,scale=0.3]{};
					\node (4g) at (7.2,-.2)  [circle,fill=black,scale=0.3]{};
					\node (4h) at (7.5,.1)  [circle,fill=black,scale=0.3]{};
					\node (4i) at (7.5,.7)  [circle,fill=black,scale=0.3]{};

					\draw [thick] (4a) -- (4b) -- (4d) -- (4a) ;
					\draw [thick] (4a) -- (4e) -- (4f) -- (4a) ;
					\draw [thick] (4a) -- (4g) -- (4h) -- (4a) ;
					\draw [thick] (4a) -- (4c) -- (4i) -- (4a) ;

					\end{tikzpicture}
					
					\caption{Three friendship graphs}
					
				\end{center}
				
			\end{figure}
			
			As a generalization of this, a \textit{multicone graph} is the join of a complete graph and multiple copies of a regular graph $H$: $K_r \triangledown sH$. Usually the graph $H$ is taken to be another complete graph, and the only multicone graphs that we consider in this paper are those of the form $K_r \triangledown sK_t$.  The friendship graph $F_s$ is thus the multicone $ K_1 \triangledown sK_2$; another example of a multicone is shown in Fig. 2.

			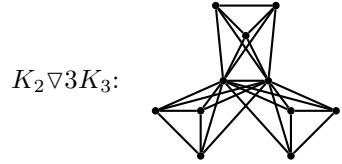
\begin{figure}
				\begin{center}
					\begin{tikzpicture}
					
					\node (C) at (0,1) [scale=1]{$K_2 \triangledown 3K_3$:};
					
					\node (1a) at (2.1,1)  [circle,fill=black,scale=0.3]{};
					\node (1b) at (2.7,1)  [circle,fill=black,scale=0.3]{};
					
					\node (2a) at (2.4,1.6)  [circle,fill=black,scale=0.3]{};
					\node (2b) at (2,2)  [circle,fill=black,scale=0.3]{};
					\node (2c) at (2.8,2)  [circle,fill=black,scale=0.3]{};
					
					\node (3a) at (1.2,.6)  [circle,fill=black,scale=0.3]{};
					\node (3b) at (1.8,.6)  [circle,fill=black,scale=0.3]{};
					\node (3c) at (1.8,0)  [circle,fill=black,scale=0.3]{};

					\node (4a) at (3.6,.6)  [circle,fill=black,scale=0.3]{};
					\node (4b) at (3,.6)  [circle,fill=black,scale=0.3]{};
					\node (4c) at (3,0)  [circle,fill=black,scale=0.3]{};
					
					\draw [thick] (1a) -- (1b);
					\draw [thick] (2a) -- (2b) -- (2c) -- (2a) ;
					\draw [thick] (3a) -- (3b) -- (3c) -- (3a) ;
					\draw [thick] (4a) -- (4b) -- (4c) -- (4a) ;
					
					\draw [thick] (1a) -- (2a);
					\draw [thick] (1a) -- (2b);
					\draw [thick] (1a) -- (2c);
					\draw [thick] (1a) -- (3a);
					\draw [thick] (1a) -- (3b);
					\draw [thick] (1a) -- (3c);
					\draw [thick] (1a) -- (4a);
					\draw [thick] (1a) -- (4b);
					\draw [thick] (1a) -- (4c);
					
					\draw [thick] (1b) -- (2a);
					\draw [thick] (1b) -- (2b);
					\draw [thick] (1b) -- (2c);
					\draw [thick] (1b) -- (3a);
					\draw [thick] (1b) -- (3b);
					\draw [thick] (1b) -- (3c);
					\draw [thick] (1b) -- (4a);
					\draw [thick] (1b) -- (4b);
					\draw [thick] (1b) -- (4c);		
					
					\end{tikzpicture}
					
					\caption{A multicone graph}
					
				\end{center}
				
			\end{figure}
	It were conjectured (see Wang, et al. \cite{WBHB, WZH}) that friendship graphs are $DAS$. Das \cite{Das} claimed to have proved this, but Abdollahi, Janbaz and Oboudi \cite{AJO} found an error in the proof, and furthermore, they proved the result for some special cases. Recently, Cioab\"a, Haemers, Vermette, and Wang \cite{CHVW} proved the conjecture for $ s \neq 16$; that is,  if $ G $ is adjacency $A$-cospectral with $ F_s$ ($ s \neq 16$), then $ G\cong F_s$.  For further information about some multicone graphs which have been characterized so far see \cite{AM, AA, AAA, AAAA, AAAAA, AAAAAA, AAAAAA1, AJO}.\\

This paper is organized as follows. In Section \ref{2}, we review some basic information and preliminaries. Then in Section \ref{3}, we state some algebraic properties about multicone graphs of the form $ K_r \triangledown sK_t $, while in Section \ref{4}, we show that these graphs are determined by their adjacency spectrum.  In Section \ref{5}, we prove that their complements are also $DAS$, and in Section \ref{6}, we prove that these graphs are also determined by their Laplacian spectrum. Finally, in Section \ref{7}, we summarize our results and propose one conjecture for further research.

			\section{Preliminary results}\label{2}
			
			\noindent	In this section, we give some results from the literature that play important roles in the rest of the paper.
			
			From both the adjacency spectrum and the Laplacian spectrum of a graph, one can deduce the number of vertices and the number of edges.  The two spectra also give additional information \cite{AJO}. We defer a similar result for the Laplacian spectrum to Section \ref{6}, where that spectrum determination  is developed.
			
			\begin{theorem}\label{the 2-1} Given a graph $G$, the following can be deduced from its adjacency spectrum:\\
				
				{\rm{(a)}} the number of closed walks of each length;
				
				{\rm{(b)}} whether or not $G$ is bipartite;
				
				{\rm{(c)}} whether or not $G$ is regular, and if so, the degree of regularity.
				
			\end{theorem}
			

				
				
				
				

			The next several results concern degrees and eigenvalues in graphs. Recall that $\Delta(G)$ (sometimes just $\Delta$) denotes the maximum degree of a vertex of a graph $G$, and similarly $\delta(G)$ (or just $\delta$) denotes the minimum degree.  If the two are different and they are the only degrees in $G$ and $\delta(G)$ is positive, then $G$ is said to be  \textit{bi-regular} or \textit{bi-degreed}.  Also, the largest eigenvalue of $G$ is called the \textit{spectral radius} (sometimes called the \textit{spectral index}) and is denoted $\rho(G)$ (or just $\rho$).

The following result, \cite{AM, AA, AAA, AAAA, HSF,  WZH} gives a bound on the spectral radius. For further information about this inequality we refer the reader to 	$ \cite{WZH} $ (see the first paragraph after Corollary 2.2 and also Theorem \ref{the 2-1} of \cite{WZH}). 
			
\begin{theorem}\label{the 2-2} If $G$ is a graph with $n$ vertices, $m$ edges, minimum degree $\delta$, and spectral radius $\rho$, then $$\rho \leq \frac{\delta - 1}{2} + \sqrt{2m - n\delta + \frac{(\delta + 1)^2}{4}}.$$
Equality holds if and only if $G$ is either regular or is bi-regular with $\Delta = n-1$.				
\end{theorem}

The next theorem gives a characterization of some graphs with three distinct eigenvalues(\cite{AM, WH}).
						
\begin{theorem}\cite{WH} \label{the 2-3}A graph has exactly one positive eigenvalue if and only if it is a complete multipartite graph with possibly some isolated vertices.	
\end{theorem}
 
The next two theorems concern regular graphs; the first can be found in Knauer \cite{Kn} and the second in Bapat \cite{Bap}.

\begin{theorem}\label{the 2-4} Let $G$ be a graph with spectral radius $\rho$.  Then the following statements are equivalent:\\
{\rm{(1)}} $G$ is regular.\\
{\rm{(2)}} $\rho$ is the average vertex degree in $G$.\\
{\rm{(3)}} $(1, 1, \ldots, 1)^{\mathsf{T}}$ is an eigenvector for $\rho$.							
\end{theorem}
						
\begin{theorem}\label{the 2-5}If $ G $ is an $r$-regular graph with eigenvalues $\lambda_1 (= r), \lambda_2, ..., \lambda_n $, then $ n-1-\lambda_1, -1-\lambda_2, ..., -1-\lambda_n $ are the eigenvalues of the complement $ \overline{G} $ of $G$.						
\end{theorem}
				
We turn now to a theorem on graphs that are not regular.

\begin{theorem}\cite{WH} \label{the 2-6} If $G$ is not regular and has exactly three eigenvalues $\theta_1 > \theta_2 > \theta_3$, then:\\
				
\noindent {\rm{(a)}} $G$ has diameter $2$;\\
{\rm{(b)}} if $\theta_1$ is not an integer, then $G$ is complete bipartite;\\
{\rm{(c)}} $\theta_2 \ge 0$ with equality if and only if $G$ is complete bipartite;\\
{\rm{(d)}} $\theta_3 < -2$.

\end{theorem}			
						
The next theorem gives the characteristic polynomial of the join of two regular graphs in terms of their individual polynomials (see also \cite{AM, AA, AAA, AAAA, CRS}).
						
\begin{theorem}\cite{CRS}\label{the 2-7} For $i = 1, 2$, let $G_i$ be an $r_i$-regular graph of order $n_i$.  Then the characteristic polynomial of their join is $$P_{G_1 \triangledown G_2}(x) = P_{G_1}(x)P_{G_2}(x) (1 - \frac{n_1 n_2}{(x-r_1)(x-r_2)}).$$	

\end{theorem}
		
The next several theorems are also on the characteristic polynomial of a graph; the first can be found in  \cite{AM, Kn}).

\begin{theorem}\label{the 2-8}

The following statements are equivalent for a nontrivial graph $G$ with characteristic polynomial $P_G(x) = \sum_{i=0}^{n}c_ix^i$, spectrum $\lambda_1 \ge \lambda_2 \ge ... \ge \lambda_n$, and spectral radius $\rho$. \\ 
{\rm{(1)}} $G$ is bipartite.\\
{\rm{(2)}} The coefficients $c_i$ for $i$ odd are all $0$.\\
{\rm{(3)}} For each $i$, $\lambda_{n+1-i} = -\lambda_i$.\\
{\rm{(4)}} $\rho = -\lambda_{n}$.

\end{theorem}
	
We note that statement $(3)$ in this theorem implies that each eigenvalue has the same multiplicity as its negative.

	The next result, including a discussion of main angles (For further information about main angles see \cite{R}), may be found in \cite{AM, AA, AAA, AAAA, CRS}).
\begin{theorem}\label{the 2-10} If $j$ is a vertex of graph $G$, then $ P_{G-j}(x)=P_{G}(x)\sum\limits_{i = 1}^m {\frac{{\alpha^2 _{ij}}}{{x - {\mu _i}}}} $, where $ m $ and $\alpha _{ij} $ are the number of distinct eigenvalues  and the main angle of graph $ G $, respectively. 
\end{theorem}

\begin{proposition}\cite{VAH} \label{prop 2-1} Let $ G $ be a disconnected graph that is determined by the Laplacian spectrum. Then the cone over $ G $, the graph $ H $; that is, obtained from $ G $ by adding one vertex that is adjacent to all vertices of $ G $, is also determined by its Laplacian spectrum. 
\end{proposition}
\section{Connected graphs $A$-cospectral with a multicone graph $ K_r\bigtriangledown sK_t $ }\label{3} 
		
In this section, we give some results on graphs that are cospectral with a multicone graph $ K_r \triangledown sK_t$. Note that the order of $K_r \triangledown sK_t$ is $r + st$, which we denote by $n$. In giving the spectrum of a graph, we often use the common notation of $[c]^k$ for an eigenvalue $c$ of multiplicity $k \ge 1$.
		
		\begin{proposition}\label{prop 3-1} If $ G $ is a graph $A$-cospectral with  multicone graph $ K_r \triangledown sK_t$, then
			
			$$\operatorname{Spec}_{A}(G)=\left\{ [-1]^{r - 1 + s(t - 1)},~[t - 1]^{s - 1}, [\frac{a + \sqrt{a^2 - 4b}}{2}]^1, [\frac{a - \sqrt{a^2 - 4b}}{2}]^1 \right\},$$
			where $a = r + t - 2 $ and $b = (r-1)(t-1)- rst. $
	
\end{proposition}		
		
\noindent \textbf{Proof} We know that $\operatorname{Spec}_{A}(K_t)=\{ [-1]^{t-1},~[t-1]^1 \}$ (see \cite{Bap}). Now,  by Theorem \ref{the 2-7} the proof is clear. $\Box$ \\

\subsection{Adjacency spectrum determination of the connected multicone graphs $K_r \triangledown sK_t$ } \label{4}

\noindent The aim of this section is to show that multicone graphs $K_r \triangledown sK_t$ are $DAS$.

\begin{lemma}\label{lem 4-1} If $ G $ is a connected graph $A$-cospectral with a multicone graph $K_r \triangledown sK_t$, then $ \delta(G) = r+t-1 $.
\end{lemma}

\noindent \textbf{Proof} Let $x = \delta(G) -(r+t-1)$. It follows from Theorem \ref{the 2-4}  that:\\

 $G$ is a regular graph if and only if $s=1$ if and only if $G$ is a complete graph.\\

 Consider the following two cases:\\

{\bf Case 1.} $s=1$. In this case $\delta(G)=r+t-1$ and there is nothing to prove.\\

{\bf Case 2.} $s\geq 2$ ($s\neq 1$). We show that $x=0$.\\

Suppose not and so $ x\neq 0 $ (in this case $\delta(G)\neq r+t-1$). It follows from Theorem \ref{the 2-2}  and Proposition \ref{prop 3-1}  that\\
\medskip 

$ \rho(G)= \frac{r+t-2+ \sqrt{8m-4n(r+t-1)+(r+t)^{2}}}{2} < \frac{r+t-2+x+ \sqrt{8m-4n(r+t-1)+(r+t)^{2}+x^{2}+(2r+2t-4n)x }}{2},$\\

\noindent where (as usual) $ n $ and $ m $ denote the numbers of vertices and edges in $ G $, respectively.  

 For convenience, we let $B = 8m-4n(r+t-1)+(r+t)^{2}$ and $C = r + t - 2n$, and also let $g(x) = x^2 + 2(r + t - 2n)x = x^2 + 2Cx.$ \\ 

Then clearly

\begin{equation*} \sqrt{B}-\sqrt{B+g(x)} < x. \tag{1} \end{equation*}

We consider the following two subcases (we show that none of the following two subcases can happen):\\

Subcase 2.1. $x < 0$.\\

Then

   $ |\sqrt{B}-\sqrt{B+g(x)}| >|x| $, since $ x<0 $.\\
   
\noindent Transposing and squaring yields

$$2B+g(x)-2\sqrt{B(B+g(x))} > x^{2}.$$

\noindent Replacing $g(x)$ by $x^2 + 2Cx$, we get 

\begin{equation*}B+Cx > \sqrt{B(B+x^2 + 2Cx)}\tag{2}.\end{equation*}

\noindent Obviously $ Cx> 0 $, since $C=r+t-2n=r+t-2(r+st)=-r+t(1-2s)<0$ and $x<0$. Squaring again and simplifying yields 

\begin{equation*} C^{2} > B. \tag{3} \end{equation*}

Therefore,

\begin{equation*} m<\frac{n(n-1)}{2}\tag{4}. \end{equation*}

Therefore, if $x<0$, then $G$ is not a complete graph. Or if $\delta(G) < r+t-1$, then $G$ is not a complete graph $(\dagger)$. On the other hand, if $x<0$ for any non-complete graph $G$ we always have $\delta(G) < r+t-1$ $(\ddagger)$. Combining $(\dagger)$ and $(\ddagger)$ we get: $\delta(G) < r+t-1$ if and only if $G$ is not a complete graph. To put that another way,  $x>0$ if and only if $G$ is  a complete graph, a contradiction, since if $G$ is a complete graph, then $x=0$.

Subcase 2.2. $x > 0$. In this case  if $G$ is non-complete graph, then  $\delta(G) > r+t-1$ (*).

On the other hand by a similar argument of Subcase 2.1 for  $x > 0$, if $\delta(G) > r+t-1$, then $G$ is not a complete graph (**). Combining (*) and (**) we have: $x<0$ if and only if $G$ is  a complete graph, a contradiction. So, we must have $x=0$. Therefore, the assertion holds. $\Box$
\begin{lemma}\label{lem 4-2} If  $G$ is a connected graph $A$-cospectral with a multicone graph $K_r \triangledown sK_t$, then it is either regular or bi-degreed with degrees $ \delta=r+t-1 $ and $\Delta=r + st - 1$.
\end{lemma}

\noindent \textbf{Proof} The result follows from Lemma \ref{lem 4-1} and Theorem \ref{the 2-2}.  $\Box$\\

In the following, we show that any connected graph $A$-cospectral with the  multicone graph $K_1 \triangledown sK_t$ is $DAS$.

\begin{lemma}\label{lem 4-3}
	If $ G $ is a connected graph $A$-cospectral with the multicone graph $K_1 \triangledown sK_t$, then $ G $ is $DAS$.
\end{lemma}

\noindent \textbf{Proof}
If $ s=1 $, there is nothing to prove, since graph $ G $ in this case is a complete graph (see  Theorem \ref{the 2-4}). Hence we suppose that $ s \neq 1 $. In this case, $ G $ is bi-degreed (see  Lemma \ref{lem 4-2}). By Lemma \ref{lem 4-2} any vertex of $ G $ is either of degree $1+t-1=t$ or $ 1 + st - 1=st$. Let $ G $ has $ \alpha $ vertices (vertex) of degree $ st $. Therefore, by  Theorem \ref{the 2-1} $ (iii) $ (sum of vertices degree of $ G $ that is sum of squares of the eigenvalues of $ G $) and Proposition \ref{prop 3-1} we have:

\begin{flushleft}
$ (\alpha)st+(st+1-\alpha)t=s(t-1)((-1)^2)+(s-1)(t-1)^2+({\frac{{ t-1 + \sqrt {(t-1)^2+4st
} }}{2}})^2+({\frac{{ t - 1 - \sqrt {(t-1)^2+4st} }}{2}})^2=st+st(t)=st(t+1)$.
\end{flushleft}
By solving the equation we get $\alpha=1 $. This means that $ G $ has one vertex of degree
$ st $, say $ j $ and $ st $ vertices of degree $ t $.   It follows from Theorem \ref{the 2-10} that
$$P_{G-j}(x)=(x-\mu_3)^{s(t-1)-1}(x-\mu_4)^{s-2}(\alpha^2_{1j}A_1+\alpha^2_{2j}A_2+\alpha^2_{3j}A_3+\alpha^2_{4j}A_4),$$ where

\begin{center}
	$ A_1=(x-\mu_2)(x-\mu_3)(x-\mu_4) $,
	
	$ A_2=(x-\mu_1)(x-\mu_3)(x-\mu_4) $,
	
	$ A_3=(x-\mu_1)(x-\mu_2)(x-\mu_4) $,
	
	$ A_4=(x-\mu_1)(x-\mu_2)(x-\mu_3) $,
	
\end{center}

\noindent with $ \mu_1=\frac{{t-1+ \sqrt {(t-1)^2+4st }}}{2} $, $\mu_2=\frac{{t-1 - \sqrt {(t-1)^2+4st }}}{2}$, $ \mu_3=-1 $, and $ \mu_4=t-1 $.
\\

As stated at the beginning of this lemma $ G $ has one vertex of degree
$ st $ and $ st $ vertices of degree $ t$. This means that graph $ G-j $ has $ st $ vertices of degree $ t-1 $. In other words, $ G-j $ is a $ (t-1) $-regular graph and it has $ st $ eigenvalues (vertices).  It is clear that  by removing the vertex $ j $ the number of edges that are deleted from graph $ G $ is  $ st=|V(G-j)| $. On the other hand, the number of the closed walks of length 2 belonging to $ G $ is:

$s(t-1)((-1)^2)+(s-1)(t-1)^2+({\frac{{ t - 1 + \sqrt {(t-1)^2+4(st)} }}{2}})^2+({\frac{{t - 1 - \sqrt {(t-1)^2+4(st)} }}{2}})^2=st+st(t)=st(t+1)$.

This means that the number of the closed walks of length 2 belonging to $ G-j $ is $ st(t+1)-2|V(G-j)|=st(t+1)-2(st)=st(t-1)$. (Or one can say that since $ G-j $ is  a $ (t-1) $-regular graph  and it has $ st $ eigenvalues, so the number of the closed walks of length 2 belonging to $ G-j $ is $ st(t-1)$).

Now, by computing the number of the closed walks of length 1 (sum of all eigenvalues that is equal to zero)  and 2 belonging to $ G-j $, we have:\\

$\gamma+\zeta+t-1=-[(s-2)\mu_4+(s(t-1)-1)\mu_3] $,

$ \gamma^2+\zeta^2+(t-1)^2=st(t-1)-[(s-2)\mu^2_4+(s(t-1)-1)\mu^2_3] $,\\

where $\gamma$ and $ \zeta$ are the eigenvalues of $G-j$. The roots are $ \gamma=t-1 $ and $ \zeta=-1 $. Therefore, $ \operatorname{Spec}_{A}(G - j) = \{\left[ { - 1} \right]^{s(t - 1)}, \left[ {t - 1} \right]^s\} = \operatorname{Spec}_{A}(sK_t) $. Hence, $ G-j\cong sK_t $, and so $ G\cong K_1 \triangledown sK_t $. $\Box$\\

Until now, we have considered only graphs $A$-cospectral with the multicone graph $ K_1 \triangledown sK_t $ (windmill-like graphs with larger sails). The next theorem extends our result to the general multicone graph $ K_r \triangledown sK_t $.

\begin{theorem}\label{the 4-1}
	If $ G $ is a connected  graph $A$-cospectral with a multicone graph $K_r \triangledown sK_t$, then $ G $ is $DAS$.
\end{theorem}

\noindent \textbf{Proof} If $s=1$ the proof is clear. Take $s\geq 2$. We perform the mathematical induction on $ r $. For $ r=1 $, the proof follows from Lemma \ref{lem 4-3}. Let the claim be true for $ r $; that is, if $ {\rm{Spec}}_{A}(M)={\rm{Spec}}_{A}(K_r\bigtriangledown sK_t) $, then $ M\cong K_r\bigtriangledown sK_t$, where $ M $ is an arbitrary graph $A$-cospectral with a multicone graph $ K_r\bigtriangledown sK_t $. We show that the claim is true for $ r+1 $; that is, we show that if $ {\rm{Spec}}_{A}(G)={\rm{Spec}}_{A}(K_{r+1}\bigtriangledown sK_t) $, then $ G\cong K_{r+1}\bigtriangledown sK_t$, where $ G $ is  a graph. It is clear that $ G $ has one vertex and $ r+st $ edges more than $ M $. By a similar argument that stated at the beginning of Lemma \ref{lem 4-3} one may deduce that $ M $ has $ r $ vertices of degree $ r+st-1 $ and $ st $ vertices of degree $ r+t-1$ and also $ G $ has $ r+1$ vertices of degree $ r+st$ and $ st $ vertices of degree $ r+t$. Hence we must have $ G\cong K_1\bigtriangledown M $, since by Lemma \ref{lem 4-2} $ G $ is bi-degreed and has $ r+1$ vertices of degree $ r+st$ and $ st $ vertices of degree $ r+t$. Now, the induction hypothesis completes the proof.  $\Box$\\

It is well-known that the smallest non-isomorphic cospectral graphs are $ \Gamma_1=C_4\cup K_1 $ and $ \Gamma_2=K_{1,4} $ (see Fig. 3). Note that $ \Gamma_1 = \overline {F_2}$ (the complement of the windmill $F_2$) is not a connected graph while $ \Gamma_2 $ is connected. Thus, we see that $ \overline {F_2} $ is not $DAS$. However, Abdollahi, Janbaz and Oboudi \cite{AJO} proved that if $ n\neq 2 $, then $ \overline{F_n}=\overline{K_1 \triangledown nK_2} $ is $DAS$. A natural question is what happens with the complement of the general multicone graph $( K_{r} \triangledown sK_t)$.  We address this in the next section.

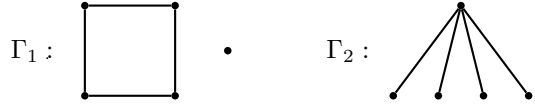
\begin{figure}
	\begin{center}
		\begin{tikzpicture}
		
		\node (A) at (-.5, .5) [circle,fill=black,scale=0.03] {$\Gamma_1$};
		\node (1) at (0,0) [circle,fill=black,scale=0.3] {};
		\node (2) at (0,1.2)  [circle,fill=black,scale=0.3] {};
		\node (3) at (1.2,1.2) [circle,fill=black,scale=0.3] {};
		\node (4) at (1.2,0)  [circle,fill=black,scale=0.3] {};
		\node (5) at (1.9,.6)  [circle,fill=black,scale=0.3] {};
		\draw [thick] (1) -- (2) -- (3) -- (4) -- (1);
		
		\node (11) at (5,1.2) [circle,fill=black,scale=0.3] {};
		\node (12) at (4.1,0)  [circle,fill=black,scale=0.3] {};
		\node (13) at (4.7,0) [circle,fill=black,scale=0.3] {};
		\node (14) at (5.3,0)  [circle,fill=black,scale=0.3] {};
		\node (15) at (5.9,0)  [circle,fill=black,scale=0.3] {};
		\draw [thick] (11) -- (12);
		\draw [thick] (11) -- (13);
		\draw [thick] (11) -- (14);
		\draw [thick] (11) -- (15);
		
		\node  at (-.7, .6)  [scale=1]{$\Gamma_1:$};
		\node  at (3.5, .6)  [scale=1]{$\Gamma_2:$};		
		
		\end{tikzpicture}
		
		\caption{A pair of $A$-cospectral graphs but non-isomorphic.}
		
	\end{center}
	
\end{figure}

\section{Graphs $A$-cospectral with complements of multicone graphs $K_r \triangledown sK_t$ }\label{5}

In this section we investigate the complements of the multicone graphs $K_r \triangledown sK_t$. Clearly, if $s = 1$, then the multicone graph is just the complete graph $K_{r+t}$, and so its complement is $(r+t)K_1$ with spectrum $\{ [0]^{r + t}\}$.  Clearly no other graph has this spectrum.  On the other hand, the case $s = 2$ is much more interesting.  The complement of $K_r \triangledown 2K_t$ is the union $rK_1 + K_{t, t}$.   The adjacency spectrum is $\{[-t], [t], [0]^{2t + r -2}\}$.  Our next theorem determines which graphs have this spectrum.

\begin{theorem} Let $G$ be a graph with adjacency spectrum  $\{[-t], [t], [0]^{2t + r -2}\}$. 
	
	${\rm{(a)}}$~~ $G$ is not connected if and only if $G \cong \overline {(K_r \triangledown 2K_t)}$. 
	
	${\rm{(b)}}$~~If $G$ is connected if and only if $G \cong K_{p, q}$, where $ p $ and $ q $  are the two roots of the  equation $x^2 -(r + 2t)x + t^2 = 0$.
\end{theorem}

\noindent \textbf{Proof}  (a) Assume that $G$ is disconnected.  Then by Theorem \ref{the 2-3}, there is a complete multipartite graph $H$ for which $G \cong H + cK_1$, where $0\leq c\leq 2t + r -2$. We show that $c=r$.   If $t\geqslant 2$, then $H$ has precisely three different eigenvalues ($L$ has two distinct eigenvalus if and only if $L=dK_n$, where $d$ and $n$ are natural numbers. Also note that $\operatorname{Spec}_{A}(K_n)=\left\{ {{{\left[ {n - 1} \right]}^1},\,{{\left[ { - 1} \right]}^{n - 1}}} \right\}$). So by Theorem \ref{the 2-8}  $H$ is a bipartite graph.  Hence $H \cong K_{t,t}$, and so $c=r$. Therefore, $G=\overline {{K_r} \triangledown 2{K_t}}$. If $t=1$, then $H\cong K_{1,1}=K_2$ and   $c=r$. The converse is clear.\\

(b) Assume that $G$ is connected.  Then $G$ cannot be regular, so by Theorem \ref{the 2-6}, it must be a complete bipartite graph $K_{p, q}$, for some $p$ and $q$.  The spectrum of $K_{p, q}$ is known to be $\{[\sqrt{pq}]^{1}, [-\sqrt{pq}]^{1}, [0]^{p+q-2}\}$ (see, for example, \cite{Kn}).  This is also the spectrum of $G \cong \overline {{K_r} \triangledown 2{K_t}}$ when $p + q = 2t +r$ and $pq = t^2$, that is, when $p$ (and likewise $q$) satisfies $x^2 -(r + 2t)x + t^2 = 0$. The converse is straightforward. $\Box$ \\

As a consequence of this theorem, we have the result that the complement of a multicone graph $\overline {{K_r} \triangledown 2{K_t}}$ is \textit{not} $DAS$.  However, these are the only graphs in this family that are not, as our next theorem shows. Before presenting the theorem, we note that $\overline{{K_r}  \triangledown s{K_t}} \cong rK_1 + K_{s(t)} $, where $K_{s(t)} $ denotes the complete $s$-partite graph with each of the partite sets being of size $t$.  We also note that the spectrum of this graph is $\{[-t]^{s-1}, [0]^{s(t-1)+r}, [t(s-1)]\}$.

\begin{theorem} 
	For $s \geq 3$, graphs  $\overline{{K_r} \triangledown s{K_t}}$ are $DAS$.
	
\end{theorem}

\noindent \textbf{Proof} Let $\operatorname{Spec}_{A}(G)=\{[-t]^{s-1}, [0]^{s(t-1)+r}, [t(s-1)]\}=\operatorname{Spec}_{A}(\overline{K_r\triangledown sK_t})$. It follows from Theorem \ref{the 2-6} that if $G$ is connected, then it is a complete bipartite graph. Therefore, it follows from Theorem \ref{the 2-8} that $ s=2 $, a contradiction. Hence $ G $ is disconnected.  Now, by Theorem \ref{the 2-3}  there is a complete multipartite graph $H$ for which $G \cong H + aK_1$, where $ s(t-1)+r\geq a\geq 1 $.  We claim that $H$ must be regular. Suppose not and so, as before, by Theorem \ref{the 2-6}  $H$ must be a complete bipartite graph, and this is impossible.  Thus, $H$ must have at least three partite sets, and since it is regular, it must be $K_s(t)$, and so $G \cong rK_1 + K_{s(t)}$, establishing the result. $\Box$ \\

\section{Laplacian spectrum determination of multicone graphs $K_r \triangledown sK_t$ }\label{6}

In this section, we consider the Laplacian spectrum of multicones.  Recall that the Laplacian matrix of a graph $G$ is the matrix $\mathbf{L}(G) = \mathbf{D}(G) - \mathbf{A}(G)$, where $\mathbf{A}(G)$ is the adjacency matrix of $G$ and $\mathbf{D}(G)$ is the degree matrix, and the Laplacian spectrum of $G$, denoted ${\rm{Spec}_L(G)}$, is the spectrum of $\mathbf{L}(G)$.

We begin with some general results of \cite{AM, AA, AAA, AAAA, Mer} on Laplacian spectra.

\begin{theorem}\label{the 6-1}
	
	Let $G$ and $H$ be graphs with Laplacian spectra $\alpha_1 \ge \alpha_2 \ge \ldots \ge \alpha_n$ and $\beta_1 \ge \beta_2 \ge \ldots \ge \beta_k$, respectively.  Then\\
	{\rm{(a)}} the Laplacian spectrum of the complement $ \overline{G} $ is $n - \alpha_1, n - \alpha_2, \ldots, n - \alpha_n-1, 0$, and \\
	{\rm{(b)}}	the Laplacian spectrum of the join $G \triangledown H$ is $n + k, k + \alpha_1, k + \alpha_2, \ldots, k + \alpha_n-1, n + \beta_1, n + \beta_2, \ldots, n + \beta_{k-1}, 0$.
\end{theorem}

\begin{theorem}\label{the 6-2}
	
	The order $n$ of a graph $G$ is a Laplacian eigenvalue of $G$ if and only if $G$ is the join of two graphs.
	
\end{theorem}

The next result gives the Laplacian spectrum of multicone graphs $K_r \triangledown sK_t$.

\begin{proposition}
	
	The Laplacian spectrum of multicone ${K_r} \triangledown s{K_t}$ is:

\begin{center}
$ \{[r+st]^{r}, [r+t]^{s(t-1)}, [r]^{s-1}, [0]^1\}$.
\end{center}
\end{proposition}
\noindent \textbf{Proof} It is clear that ${{\rm{Spec}}_L(K_t) = \{[t]^{t-1}, [0]^1\}}$ and so ${{\rm{Spec}}_L(sK_t) = \{[t]^{s(t-1)}, [0]^s\}}$. Now, By Theorem \ref{the 6-1} (b) the proof is straightforwad.$\Box$
\begin{theorem}
	
 Multicone graphs ${K_r} \triangledown s{K_t}$ are $DLS$.
	
\end{theorem}

\noindent \textbf{Proof} If $s=1$, the proof is clear. So, we consider $s\geq 2$. The proof is by induction on $r$. By Proposition \ref{prop 2-1} the result is clearly true when $r=1$.  Assume that the theorem holds for $r$; that is, if ${\rm{Spec}}_{L}(G)={\rm{Spec}}_{L}({K_r} \triangledown s{K_t})=\{[r+st]^r, [r+t]^{s(t-1)}, [r]^{s-1}, [0]^1\}$, then $G \cong {K_r} \triangledown s{K_t}.$ We show that if  ${\rm{Spec}}_{L}(H)={\rm{Spec}}_{L}({K_{r+1} \triangledown s{K_t})=\{[r+st+1]^{r+1}, [r+t+1]^{s(t-1)}, [r+1]^{s-1}, [0]^1\}}$, then $H \cong {K_{r+1}} \triangledown s{K_t}.$  It follows from Theorem \ref{the 6-2} that $H$ and $G$ are the join of two graphs.  On the other hand, $H$ has one vertex, say $e$ and $r+st$ edges more  than $G$. By Theorem \ref{the 6-1} (a) ${\rm{Spec}}_{L}(\overline{H})=\{[0]^{r+2}, [st-t]^{s(t-1)}, [st]^{s-1}\}$. Therefore, $\overline{H}$ has $r+2$ connected components, since it has $r+2$ eigenvalues that are zero. Note that if one can receive to $ H $ from $ G $, then this means that $G\subseteq H$ and conversely.\\

 We prove that  only if $e$ join to $G$ one can receive to $H$ from $G$. Suppose not and so let one can receive to $ H $ from $ G $ without joining $ e $ to $ G $. Since $ e $ does not join to $ G={K_r} \triangledown s{K_t}$, so, there is  a vertex of $ G $ that is not adjacent to $ e $, say $o$. In this case connected components of $ \overline{H} $ are as following:\\
1. Let $o$ be belonging to $K_r$. It is clear that two vertices of two components of $sK_t$ must be adjacent (because of the number of edges of $H$), say $uv$. Then  connected components of $ \overline{H} $ are $ K_{1, 1} $, $K_{s(t)}-uv$ and $(r-1)K_1 $. We prove that there is no graph $\overline{H}$ with these connected components. Suppose not and first let $K_{s(t)}=K_2=\overline{2K_1}$. Then connected components of $ \overline{H} $ are $ K_{1, 1} $ and $(r+1)K_1$.  Therefore, $\overline{H}=(r+1)K_1+K_2$ or $H=2K_1\bigtriangledown K_{r+1}$, a contradiction, since $e$ is not adjacent with all vertices of $K_r$. For another cases $ \overline{H} $ has $r+1$ connected components, a contradiction.\\

2. Let $o$ be belonging to $sK_t$.  Then  connected components of $ \overline{H} $ are $ rK_1 $ and $K_{s(t)}-uv\ast e$, where $K_{s(t)}-uv\ast e$ means that $ e$ is adjacent to one of  vertices $K_{s(t)}-uv$. If $K_{s(t)}=K_2=\overline{2K_1}$, then connected components of $ \overline{H} $ are $ K_{1, 1} $ and $(r+1)K_1$ and so $\overline{H}=(r+1)K_1+K_2$ or $H=K_{r+1}\bigtriangledown2K_1$,   a contradiction, since $e$ is not adjacent with all vertices of $2K_1$. For another cases $ \overline{H} $ has $r+1$ connected components, a contradiction.\\
Hence one may deduce that:\\

$G$ is a subgraph of $H$ if and only if $H$ can  be obtained from $G$ if and only if $e$ join to $G$.\vspace{1mm}

 So, $H$ can only be obtained from $G$ by joining $e$ to $G$.  In other words,  $ H=K_1\bigtriangledown G $, and then by the induction hypothesis the proof is straightforward.$\Box$\\

Since friendship graphs form a special family of multicone graphs, we have the following result.

\begin{corollary}
	
The friendship graph $F_s={K_1} \triangledown s{K_2}$ is $DLS$.
\end{corollary}

\section{Conclusion remarks and a Conjecture}\label{7}

The following theorem summarizes the most important results in this paper regarding the adjacency and Laplacian spectrum determination of multicone graphs.

\begin{theorem}
	
	\smallskip
	
	${\rm{(a)}}$~~For all $r, s,$ and $t$, the connected multicone graph ${K_r} \triangledown s{K_t}$ is both  $DAS$ and $DLS$.
	
	${\rm{(b)}}$~~For all $r, s,$ and $t$ with $s>2$, the graphs $\overline {{K_r} \triangledown s{K_t}}$ are both  $DAS$ and $DLS$.
\end{theorem}

Of course, friendship graphs are an interesting, and in some ways exceptional, special type of multicone, consisting as they do, of a collection of triangles with one common vertex.  Cioab\"a, et al. \cite{CHVW}
proved that the friendship graph $F_{16}$ is not $DAS$, and, as we observed earlier, neither is the complement of $F_2$, and these are the only exceptions.

\begin{remark}

	                    ${\rm{(a)}}$~~For $s \neq 16$, the friendship graph $F_s=K_1 \triangledown sK_2$ is $DAS$ and for any $s$, $F_s$ is $DLS$.
	
	                   ~~~~~~~~~ ${\rm{(b)}}$~~For $s \neq 2$, the friendship graph $\overline{F_s=K_1 \triangledown sK_2}$ is   $DAS$ and for any $s$, $\overline{F_s}$ is $DLS$.

                     ~~~~~~~~~${\rm{(c)}}$~~ Consider $ K_{1, 3}\bigtriangledown K_{r-1}=3K_1\bigtriangledown K_r$  and $(K_3\cup K_1)\bigtriangledown K_{r-1} $ that have the same signless Laplacian spectrum (see Corollary 2.2 of \cite{LP} and Theorem 2.1 of \cite{LP}) but are non-isomorphic.\\

\end{remark}

As we noted at the beginning of the paper, a third type of matrix that gives the adjacencies of a graph has been studied, the signless Laplacian matrix, defined as ${\bf{S}}(G) = {\bf{D}}(G) + {\bf{A}}(G)$ (in contrast to the ordinary Laplacian ${\bf{S}}(G) = {\bf{D}}(G) - {\bf{A}}(G)$), with the corresponding \textit{determined by signless Laplacian spectrum}.  Friendship graphs are known to be $DQS$.  Thus, we conclude with the following conjecture.

\begin{conjecture}   Multicone graphs $K_r \triangledown sK_t $,  except in multicone graphs $K_r\bigtriangledown 3K_1$,  are $DQS$.
\end{conjecture}

\end{document}